\documentclass[pdftex]{amsart}
\usepackage{amsmath, amssymb}
\usepackage{type1cm}
\usepackage{tikz, tikz-cd}
\usetikzlibrary{matrix,arrows,shapes,calc}

\newtheorem{defi}{Definition}[section]
\newtheorem{thm}[defi]{Theorem}
\newtheorem{prop}[defi]{Proposition}
\newtheorem{lem}[defi]{Lemma}
\newtheorem{cor}[defi]{Corollary}
\newtheorem{eg}[defi]{Example}
\newtheorem{rem}[defi]{Remark}
\newtheorem{claim}[defi]{Claim}

\newtheorem{Q}[defi]{Question}
\newtheorem{pf}{Proof}

\DeclareMathOperator{\D}{D^{\mathrm{b}}}

\DeclareMathOperator{\Ext}{Ext}

\DeclareMathOperator{\Ho}{H}
\DeclareMathOperator{\homo}{Hom}

\DeclareMathOperator{\map}{\longrightarrow}

\DeclareMathOperator{\R}{R}
\DeclareMathOperator{\RHom}{RHom}

\DeclareMathOperator{\Sym}{Sym}
\DeclareMathOperator{\Spec}{Spec}
\newcommand{\stsh}{\mathcal{O}}

\title[exceptional collections on toric varieties]{strong full exceptional collections on certain toric varieties with picard number three via mutations}
\author{wahei hara}
\address{Department of Mathematics, School of Science and Engineering, Waseda University, 3-4-1 Ohkubo, Shinjuku, Tokyo 169-8555, Japan}
\email{waheyhey@ruri.waseda.jp}
\date{}

\begin{document}

\begin{abstract}
In this paper, we study derived categories of certain toric varieties with Picard number three that are blowing-up another toric varieties along their torus invariant loci of codimension at most three. We construct strong full exceptional collections by using Orlov's blow-up formula and mutations.
\end{abstract}

\maketitle

\tableofcontents

\section{Introduction}
An object $E$ of a triangulated category $\mathcal{D}$ is called \textit{exceptional} if
\begin{equation*}
\homo_{\mathcal{D}}(\mathcal{E}, \mathcal{E}[i]) = \left\{
\begin{array}{ll}
\Bbbk  & \mbox{if $i = 0$,} \\
0  & \mbox{if $i \neq 0$,}
\end{array}
\right.
\end{equation*}
and a sequence of exceptional objects $E_1, \dots, E_r$ is called \textit{full exceptional collection} if they generate whole category $\mathcal{D}$ and $\homo_{\mathcal{D}}(\mathcal{E}_l, \mathcal{E}_k [i]) = 0$ for all $1 \leq  k < l \leq r$ and all $i \in \mathbb{Z}$. In addition, the full exceptional collection is \textit{strong} if $\homo_{\mathcal{D}}(\mathcal{E}_k, \mathcal{E}_l[i]) = 0$ for all $1 \leq k < l \leq r$ and all $i \neq 0$.

If one finds a full exceptional collection, one can draw many information on the triangulated category $\mathcal{D}$. However, a triangulated category does not always admit a full exceptional collection. For example, derived categories of Calabi-Yau  varieties do not have any exceptional collection. For toric projective case, Y. Kawamata proved in \cite{Kawamata06} that:

\begin{thm}[\cite{Kawamata06}]
For any smooth projective toric Deligne-Mumford stack $\mathcal{X}$, its derived category $\D(\mathcal{X})$ has a full exceptional collection.
\end{thm}

About the existence of \textit{strong} full exceptional collections, there are same conjectures. The following question is due to A. King \cite{King97}.

\begin{Q} \label{question}
For any smooth toric variety $X$, does its derived category $\D(X)$ have a strong full exceptional collection consisting of line bundles?
\end{Q}

However,  the answer of this question is negative in general:
\begin{enumerate}

\item[(1)] First, L. Hille and M. Perling constructed in \cite{HP06} a $2$-dimensional counterexample to Question \ref{question}. 
This example is just the Hirzebruch surface $\mathbb{F}_2$ iteratively blown-up three times.

\item[(2)] Further, M. Macha{\l}ek presented an infinite list of counterexamples for Question \ref{question} in \cite{M11}.

\item[(3)] A. Efimov showed in \cite{Efimov14} that there are infinitely many counterexamples for Question \ref{question} 
that are smooth toric \textit{Fano} varieties with \textit{Picard number three}.

\end{enumerate}

Toric varieties with Picard number at most two are studied by L. Costa, R.M. Mir\'{o}-Roig \cite{CM04}, and  they proved that their derived categories have strong full exceptional collections consisting of line bundles. A. Day, M. Laso\'{n}, M. Micha{\l}ek \cite{DLM09}, L. Costa, R.M. Mir\'{o}-Roig \cite{CM12}, and M. Laso\'{n}, M. Micha{\l}ek \cite{LM11} studied the derived categories of toric varieties with Picard number three that are blowing-up of another toric varieties along codimension two loci. In this paper, we generalize their results and newly study the toric varieties which are blowing-up of another toric varieties along codimension three loci. More precisely, we prove the following theorem.

\begin{thm}[= \ref{Main Theorem}]
Let $X$ be a smooth projective toric variety with Picard number two, and $\widetilde{X}$ a blowing-up of $X$ along a torus invariant closed subvariety $Y \subset X$. If the codimension of $Y$ in $X$ is at most three, then $\D(\widetilde{X})$ has a strong full exceptional collection consisting of line bundles.
\end{thm}

In the previous works \cite{DLM09, CM12, LM11}, the authors used Bondal's \textit{Frobenius splitting method} to construct a strong full exceptional collection consisting of line bundles in a derived category of a toric variety.  In this paper, we take different approach, namely we prove the theorem by Orlov's blow-up formula and the \textit{mutation method}. If we use the Frobenius splitting method, we need to check that the collection as an output is actually full, exceptional, and strong. But in our case, because the operation of \textit{mutation} keeps the condition ``full exceptional", what we need to check is only the strongness of the collection. This makes the computations in the proof much easier and more elementary, and also enables us to generalize the previously known results. Note that the difficulty of the mutation method is the difficulty of the explicit calculations of mutated objects, but we find a new procedure of mutation operations which we can easily calculate.

\vspace{0.1in}

\noindent
\textbf{Acknowledgements. }
The author would like to express his gratitude to his supervisor Yasunari Nagai for beneficial conversations and helpful advices. He is also grateful to the referee for careful reading and giving useful comments.

\section{Preliminaries}

Let $\Bbbk$ be an algebraic closed field of any characteristic.

\subsection{Semiorthogonal decompositions and exceptional collections}
Let $\mathcal{D}$ be a triangulated category over a field $\Bbbk$.
\begin{defi}
\rm
Let $\mathcal{A}_1, \dots, \mathcal{A}_r$ be triangulated full subcategories of $\mathcal{D}$. The sequence of subcategories $\mathcal{A}_1, \dots, \mathcal{A}_r$ is called a \textit{semiorthogonal collection} in $\mathcal{D}$ if $\homo_{\mathcal{D}}(\mathcal{F}, \mathcal{E}) = 0$ for all $1 \leq i < j \leq r$ and all $\mathcal{E} \in \mathcal{A}_i, \mathcal{F} \in \mathcal{A}_j$. A semiorthogonal collection $\mathcal{A}_1, \dots, \mathcal{A}_r$ is called a \textit{semiorthogonal decomposition} if it generates the whole category $\mathcal{D}$, i.e.\ if the smallest triangulated subcategoy of $\mathcal{D}$ that contains all subcategories $\mathcal{A}_1, \dots, \mathcal{A}_r$ coincides with $\mathcal{D}$.  In such case, we write
\[ \mathcal{D} = \langle \mathcal{A}_1, \dots, \mathcal{A}_r \rangle. \]
\end{defi}

\begin{defi}
\rm
\begin{enumerate}
\item[(i)] An object $\mathcal{E} \in \mathcal{D}$ is called an \textit{exceptional object} if 
\begin{equation*}
\homo_{\mathcal{D}}(\mathcal{E}, \mathcal{E}[i]) = \left\{
\begin{array}{ll}
\Bbbk  & \mbox{if $i = 0$,} \\
0  & \mbox{if $i \neq 0$.}
\end{array}
\right.
\end{equation*}
\item[(ii)] A sequence of exceptional objects $\mathcal{E}_1, \dots, \mathcal{E}_r$ is called an \textit{exceptional collection} if $\homo_{\mathcal{D}}(\mathcal{E}_l, \mathcal{E}_k [i]) = 0$ for all $1 \leq k < l \leq r$ and all $i \in \mathbb{Z}$.
\item[(iii)] An exceptional collection $\mathcal{E}_1, \dots, \mathcal{E}_r$ is \textit{full} if it generates the whole category $\mathcal{D}$. In such case, we write
\[ \mathcal{D} = \langle \mathcal{E}_1, \dots, \mathcal{E}_r \rangle. \]
\item[(iv)] An exceptional collection $\mathcal{E}_1, \dots, \mathcal{E}_r$ is \textit{strong} if $\homo_{\mathcal{D}}(\mathcal{E}_k, \mathcal{E}_l[i]) = 0$ for all $1 \leq k < l \leq r$ and all $i \neq 0$.
\end{enumerate}
\end{defi}

\begin{eg}[\cite{Beilinson79}]
\rm
An $n$-dimensional projective space $\mathbb{P}^n$ has a strong full exceptional collection consisting of line bundles called Beilinson collection
\[  \D(\mathbb{P}^n) = \langle \stsh, \stsh(1), \stsh(2), \dots, \stsh(n) \rangle. \]
\end{eg}

\begin{rem} \label{SOD FEC}
\rm
If $\mathcal{E} \in \mathcal{D}$ is an exceptional object, the category $\langle \mathcal{E} \rangle$ generated by $\mathcal{E}$ is equivalent to the derived category of a point $\D(\Bbbk) = \D(\Spec \Bbbk)$. If a sequence of objects $\mathcal{E}_1, \dots, \mathcal{E}_r$ is a full exceptional collection in $\mathcal{D}$, then a sequence of subcategories $\langle \mathcal{E}_1 \rangle, \dots, \langle \mathcal{E}_r \rangle$ is a semiorthogonal decomposition of  $\mathcal{D}$. Conversely, if the sequence of subcategories $\mathcal{A}_1, \dots, \mathcal{A}_r$ is a semiorthogonal decomposition of $\mathcal{D}$ and each subcategory $\mathcal{A}_i$ has a full exceptional collection, then $\mathcal{D}$ also has a full exceptional collection.
\end{rem}

\begin{rem}
\rm
If an $\mathrm{Ext}$-finite category $\mathcal{D}$ (which means that for any $\mathcal{F}, \mathcal{F}' \in \mathcal{D}$ the vector space $\oplus_{i \in \mathbb{Z}} \homo_{\mathcal{D}}(\mathcal{F}, \mathcal{F}'[i])$ is finite dimensional) has a strong full exceptional collection $\mathcal{E}_1, \dots, \mathcal{E}_r$, then  there is an equivalence from $\mathcal{D}$ to the derived category of right modules over the non-commutative ring $A = \mathrm{End}\left( \bigoplus_{i=1}^r \mathcal{E}_i \right)$ defined by 
\[ \RHom(\bigoplus \mathcal{E}_i, - ) : \mathcal{D} \map \D(\text{$\mathrm{mod}$-}A). \]
This equivalence was first proved by A. Bondal  in \cite{Bondal90} when $\mathcal{D}$ is a derived category of a smooth projective variety with a strong full exceptional collection.
\end{rem}

\subsection{Mutations}

For an object $\mathcal{E} \in \mathcal{D}$, we define subcategories $\mathcal{E}^{\bot}, {}^{\bot}\mathcal{E} \subset \mathcal{D}$ by
\begin{align*}
\mathcal{E}^{\bot} &:= \{ \mathcal{F} \in \mathcal{D} \mid \homo_{\mathcal{D}}(\mathcal{E}, \mathcal{F}) = 0 \} \\
{}^{\bot}\mathcal{E} &:= \{ \mathcal{F} \in \mathcal{D} \mid \homo_{\mathcal{D}}(\mathcal{F}, \mathcal{E}) = 0 \}.
\end{align*}

\begin{defi}
\rm
Let $\mathcal{E} \in \mathcal{D}$ be an exceptional object. For an object $\mathcal{F}$ in ${}^{\bot}\mathcal{E}$, we define the \textit{left mutation of $\mathcal{F}$ through $\mathcal{E}$} as the object $\mathbb{L}_{\mathcal{E}}(\mathcal{F})$ in $\mathcal{E}^{\bot}$ that lies in an exact triangle
\[ \RHom(\mathcal{E}, \mathcal{F}) \otimes \mathcal{E} \map \mathcal{F} \map \mathbb{L}_{\mathcal{E}}(\mathcal{F}). \]
Similarly, for an object $\mathcal{G}$ in $\mathcal{E}^{\bot}$, we define the \textit{right mutation of $\mathcal{G}$ through $\mathcal{E}$} as the object $\mathbb{R}_{\mathcal{E}}(\mathcal{G})$ in ${}^{\bot}\mathcal{E}$ which lies in an exact triangle
\[ \mathbb{R}_{\mathcal{E}}(\mathcal{G}) \map \mathcal{G} \map \RHom(\mathcal{G}, \mathcal{E})^* \otimes \mathcal{E}. \]
\end{defi}

\begin{lem}[\cite{Bondal90}]
Let $\mathcal{E}_1, \mathcal{E}_2$ be an exceptional pair (i.e.\ an exceptional collection consisting of two objects). Then, the following holds.
\begin{enumerate}
\item[(i)] The left (resp. right) mutated object $\mathbb{L}_{\mathcal{E}_1}(\mathcal{E}_2)$ (resp. $\mathbb{R}_{\mathcal{E}_2}(\mathcal{E}_1)$) is again an exceptional object.
\item[(ii)] The pairs of exceptional objects $\mathcal{E}_1, \mathbb{R}_{\mathcal{E}_1}(\mathcal{E}_2)$ and $\mathbb{L}_{\mathcal{E}_2}(\mathcal{E}_1), \mathcal{E}_2$ are again exceptional pairs.
\end{enumerate}
Let $\mathcal{E}_1, \dots, \mathcal{E}_r$ be a full exceptional collection in $\mathcal{D}$. Then
\begin{enumerate}
\item[(iii)] The collection 
\[ \mathcal{E}_1, \dots, \mathcal{E}_{i-1}, \mathbb{L}_{\mathcal{E}_i}(\mathcal{E}_{i+1}), \mathcal{E}_i, \mathcal{E}_{i+2}, \dots, \mathcal{E}_r \]
is again full exceptional for each $1 \leq i \leq r-1$. Similarly, the collection 
\[ \mathcal{E}_1, \dots, \mathcal{E}_{i-2}, \mathcal{E}_i, \mathbb{R}_{\mathcal{E}_i}(\mathcal{E}_{i-1}), \mathcal{E}_{i+1}, \dots, \mathcal{E}_r \] 
is again full exceptional for each $2 \leq i \leq r$.
\end{enumerate}
\end{lem}

\begin{lem}[\cite{Bondal90}] \label{lem mutation}
\begin{enumerate}
\item[(i)] Let $\mathcal{E}_1, \mathcal{E}_2 \in \mathcal{D}$ be an exceptional pair. Assume that we have $\homo_{\mathcal{D}}(\mathcal{E}_1, \mathcal{E}_2[i]) = 0$ for all $i \in \mathbb{Z}$. Then, $\mathbb{L}_{\mathcal{E}_1}(\mathcal{E}_2) \simeq \mathcal{E}_2$ and $\mathbb{R}_{\mathcal{E}_2}(\mathcal{E}_1) \simeq \mathcal{E}_1$.
\item[(ii)] Let $\D(X) = \langle \mathcal{E}_1, \mathcal{E}_2, \dots, \mathcal{E}_{r-1}, \mathcal{E}_r \rangle$ be an full exceptional collection in a derived category of smooth projective variety $\D(X)$. Then, the following two collections
\begin{align*}
\langle \mathcal{E}_2, \dots, \mathcal{E}_{r-1}, \mathcal{E}_r, \mathcal{E}_1 \otimes \omega_X^{-1} \rangle,
~~~ \langle \mathcal{E}_r \otimes \omega_X, \mathcal{E}_1, \mathcal{E}_2, \dots, \mathcal{E}_{r-1} \rangle
\end{align*}
are also full exceptional collections in $\D(X)$.
\end{enumerate}
\end{lem}

\subsection{Orlov's formulas}

We recall Orlov's two formulas that give semiorthogonal decompositions of derived categories.
We will use these formulas to construct a full exceptional collection on the derived category of our toric variety.

\begin{thm}[\cite{Orlov93}] \label{projective bundle}
Let $X$ be a smooth projective variety and $\mathcal{E}$ a vector bundle of rank r+1 on $X$. Consider the projectivization of $\mathcal{E}$, $p : \widetilde{X} := \mathbb{P}_X(\mathcal{E}) \to X$. Then, the functor $p^* : \D(X) \to \D(\widetilde{X})$ is fully faithful, and $\D(\widetilde{X})$ has a semiorthogonal decomposition
\[ \D(\widetilde{X}) = \langle p^*\D(X), p^*\D(X) \otimes \stsh_{p}(1), \dots, p^*\D(X) \otimes \stsh_{p}(r) \rangle \]
where $\stsh_{p}(1)$ is the tautological line bundle of $\mathbb{P}_X(\mathcal{E})$.
\end{thm}

\begin{thm}[\cite{Orlov93}] \label{blow up}
Let $X$ be a smooth projective variety, and $Y \subset X$ a smooth closed subvariety of codimension $c ~ ( \geq 2)$. Let $f : \widetilde{X} := \mathrm{Bl}_Y X\to X$ be  a blowing-up of $X$ along $Y$ and $E$ its exceptional divisor,
\[ \begin{tikzcd}
    E = \mathbb{P}(\mathcal{N}^*_{Y/X}) \arrow[hook]{r} \arrow{d}[swap]{\pi} & \widetilde{X} \arrow{d}{f} \\
    Y \arrow[hook]{r}  & X  
\end{tikzcd} \]

Then, the functors $f^* : \D(X) \to \D(\widetilde{X})$ and $\iota_*\pi^* : \D(Y) \to \D(\widetilde{X})$ are fully faithful, and $\D(\widetilde{X})$ has a semiorthogonal decomposition
\[ \D(\widetilde{X}) = \langle\iota_*\pi^*\D(Y) \otimes \stsh((c-1)E), \dots, \iota_*\pi^*\D(Y) \otimes \stsh(E), f^*\D(X) \rangle. \]
\end{thm}

\section{Main theorem and comparison with known results}

First, we recall the following result due to L. Costa and R.M. Mir\'{o}-Roig.

\begin{prop}[\cite{CM04}] \label{toric projspbdl}
Let $X$ be a smooth projective toric variety, and $\mathcal{E}$ a vector bundle of rank $r+1$ on $X$ whose projectivization $Z = \mathbb{P}_X(\mathcal{E})$ is also toric. Assume that $\D(X)$ has a full exceptional collection consisting of line bundles, then $\D(Z)$ also has a full exceptional collection consisting of line bundles. Moreover, if the full exceptional collection in $\D(X)$ is strong, then $\D(Z)$ has a strong full exceptional collection consisting of line bundles.
\end{prop}

A smooth projective toric variety with Picard number one is just a projective space.
On the other hand, the geometric structure of smooth projective toric varieties with Picard number two is given by the following theorem.

\begin{thm}[\cite{CLS}, \cite{Kleinschmidt88}] \label{toric pic2}
Let $X$ be a smooth projective toric variety with Picard number two. Then, there are integers $s, r \geq 1$, $s+r = \dim X$, and $0 \leq a_1 \leq a_2 \leq \cdots \leq a_r$ such that
\[ X \simeq \mathbb{P}_{\mathbb{P}^s}(\stsh_{\mathbb{P}^s} \oplus \stsh_{\mathbb{P}^s}(a_1) \oplus \stsh_{\mathbb{P}^s}(a_2) \oplus \cdots \oplus \stsh_{\mathbb{P}^s}(a_r)). \]
\end{thm}

From the above, we have the following.

\begin{cor}
Question \ref{question} is true for smooth toric varieties with Picard number two, i.e. their derived categories have strong full exceptional collections consisting of line bundles.
\end{cor}

For toric varieties with Picard number three, Question \ref{question} is not true in general. 
More precisely, A. Efimov proved in \cite{Efimov14} that there are infinitely many smooth toric Fano varieties with Picard number three whose derived categories do not have strong full exceptional collections consisting of line bundles.
\\
\\
Our main theorem is as follow.

\begin{thm}[Main Theorem] \label{Main Theorem}
Let $X$ be a smooth projective toric variety with Picard number two, and $\widetilde{X}$ a blowing-up of $X$ along a torus invariant closed subvariety $Y \subset X$. 
If the codimension of $Y$ in $X$ is at most three, then $\D(\widetilde{X})$ has a strong full exceptional collection consisting of line bundles.
\end{thm}

\begin{rem} \rm
There are classifications for toric Fano threefolds and toric Fano fourfolds by V. Batyrev and H. Sato \cite{Batyrev99, Sato00}. 
Using these classifications, A. Bernardi, S. Tirabassi, and H. Uehara proved that Question \ref{question} is true for all toric Fano threefolds \cite{BT09, Uehara14},
 and N. Prabhu-Naik did for all toric Fano fourfolds \cite{Prabhu-Naik15}. Their method of the proof is the Bondal's Frobenius splitting method, and the last author also used some computational tools. Our Theorem and Proposition \ref{toric projspbdl} give another proof of their results for all toric Fano threefolds with Picard number three and for 27 (i.e. all except one) toric Fano fourfolds with Picard number three, without using the Frobenius splitting method.
\end{rem}

\begin{rem} \rm
There are some previous works about the Question \ref{question} for toric variety with Picard number three \cite{DLM09, CM12, LM11}. Our theorem includes these previous results.
\end{rem}

\section{Some lemmas}

To prove the theorem, we will use  the following lemmas.

\begin{lem} \label{lemA}
Let $X$ be an $n$-dimensional smooth projective variety, and Y a smooth closed subvariety of $X$ of codimension $c ~ (\geq 2)$. 
Let $\widetilde{X} := \mathrm{Bl}_Y X$ be a blowing-up of $X$ along $Y$, $E$ the exceptional divisor, 
$\iota : E \hookrightarrow \widetilde{X}$ the closed immersion,
$f : \widetilde{X} \to X$ the projection,
and $\pi : E \to Y$ the restriction of $f$ on $E$.
If $\mathcal{L}$ and $\mathcal{M}$ are line bundles on $X$ and on $Y$, respectively,
then there is a natural isomorphism
\[ \Ext_{\widetilde{X}}^i(\iota_*\pi^*\mathcal{M} \otimes \stsh(kE), f^* \mathcal{L}) \simeq \Ext^{i-1}_Y(\mathcal{M}, \mathcal{L}|_{Y} \otimes \Sym^{k-1} \mathcal{N}_{Y/X}^*) \]
for $1 \leq k \leq c-1$, where $\mathcal{N}_{Y/X}$ is the normal bundle of $Y \subset X$.
\end{lem}

\begin{pf}
\rm
By Serre duality and the projection formula, we have
\begin{align*}
&\Ext_{\widetilde{X}}^i(\iota_*\pi^*\mathcal{M} \otimes \stsh(kE), f^* \mathcal{L}) \\
\simeq ~ &\Ext_{\widetilde{X}}^{n-i}(f^* \mathcal{L}, \iota_*\pi^*\mathcal{M} \otimes \stsh(kE) \otimes f^*\omega_X \otimes \stsh((c-1)E))^* \\
\simeq ~ &\Ho^{n- i}(\mathbb{P}_{Y}(\mathcal{N}_{Y/X}^*) , \pi^*(\mathcal{L}^*|_{Y} \otimes \mathcal{M} \otimes {\omega_X}|_{Y}) \otimes \stsh_{\pi}(-k - c + 1))^*.
\end{align*}
By using the Leray spectral sequence
\begin{align*}
E_2^{p,q} = \Ho^p&(Y, \mathcal{L}^*|_{Y} \otimes \mathcal{M} \otimes {\omega_X}|_{Y} \otimes \R^q\pi_*(\stsh_{\pi}(-k - c + 1))) \\
&\Rightarrow E^{p+q} = \Ho^{p+q}(\mathbb{P}_{Y}(\mathcal{N}_{Y/X}^*) , \pi^*(\mathcal{L}^*|_{Y} \otimes \mathcal{M} \otimes {\omega_X}|_{Y}) \otimes \stsh_{\pi}(-k - c + 1))
\end{align*}
and the formula
\begin{equation*}
\R^q\pi_*(\stsh_{\pi}(-k - c + 1)) = \left\{
\begin{array}{ll}
\pi_*(\stsh_{\pi}(k-1))^* \otimes \bigwedge^c \mathcal{N}_{Y/X} & \mbox{if $q = c-1$,} \\
0  & \mbox{otherwise}
\end{array}
\right.
\end{equation*}
(Note that $-k - c +1 \leq -c$), we obtain an isomorphism
\begin{align*}
&\Ho^{n- i}(\mathbb{P}_{Y}(\mathcal{N}_{Y/X}^*) , \pi^*(\mathcal{L}^*|_{Y} \otimes \mathcal{M} \otimes {\omega_X}|_{Y}) \otimes \stsh_{\pi}(-k - c + 1))^* \\
\simeq ~ &\Ho^{n-c-i+1}(Y, \mathcal{L}^*|_{Y} \otimes \mathcal{M} \otimes {\omega_X}|_{Y} \otimes \pi_*(\stsh_{\pi}(k-1))^* \otimes \bigwedge^c \mathcal{N}_{Y/X})^*
\end{align*}
Again, by using Serre duality and the adjunction formula $\omega_Y \simeq {\omega_X}|_{Y} \otimes \bigwedge^c \mathcal{N}_{Y/X}$, we have
\begin{align*}
&\Ho^{n-c-i+1}(Y, \mathcal{L}^*|_{Y} \otimes \mathcal{M} \otimes {\omega_X}|_{Y} \otimes \pi_*(\stsh_{\pi}(k-1))^* \otimes \bigwedge^c \mathcal{N}_{Y/X})^* \\
\simeq ~ &\Ho^{i-1}(Y, \mathcal{L}|_{Y} \otimes \mathcal{M}^* \otimes \pi_*(\stsh_{\pi}(k-1))) \\
\simeq ~ &\Ext^{i-1}(\mathcal{M}, \mathcal{L}|_{Y} \otimes \mathrm{Sym}^{k-1}\mathcal{N}_{Y/X}^*).
\end{align*}
Therefore, we obtain the desired isomorphism. \qed
\end{pf}

Recall that a line bundle $\mathcal{L}$ on $X$ is \textit{acyclic} if $\Ho^i(X, \mathcal{L}) = 0$ for all $i \neq 0$.

\begin{lem} \label{lemB}
Let $X$, $Y$, $\widetilde{X}$, and $E$ as above. If $\mathcal{L}$ is an acyclic line bundle on $X$, then the line bundle $f^* \mathcal{L} \otimes \stsh(kE)$ on $\widetilde{X}$ is acyclic for $0 \leq k \leq c-1$.
\end{lem}

\begin{pf}
\rm
When $k =0$, the claim follows from the projection formula. Let us assume that $k \geq 1$ and $f^*\mathcal{L} \otimes \stsh((k-1)E)$ is acyclic. Let us consider the fundamental sequence
\[ 0 \to f^*\mathcal{L} \otimes \stsh((k-1)E) \to f^*\mathcal{L} \otimes \stsh(kE) \to f^*\mathcal{L}|_{E} \otimes \stsh_E(kE) \to 0. \]
Since $\R\pi_*\stsh_E(kE) = R\pi_*\stsh_{\pi}(-k) = 0$ for $1 \leq k \leq c-1$,
we have $f^*\mathcal{L} \otimes \stsh(kE)$ is also acyclic. \qed
\end{pf}

\section{Proof of Theorem \ref{Main Theorem}, codimension two case}

By Theorem \ref{toric pic2}, we may assume that a toric variety $X$ of Picard number two is a projective space bundle over a projective space $\mathbb{P}^s$.
Let $\mathcal{E} = \stsh_{\mathbb{P}^s} \oplus \stsh_{\mathbb{P}^s}(a_1) \oplus \cdots \oplus \stsh_{\mathbb{P}^s}(a_r) ~ (0 \leq a_1 \leq \cdots \leq a_r)$ be a vector bundle on $\mathbb{P}^s$
such that $X = \mathbb{P}_{\mathbb{P}^s}(\mathcal{E})$.
Fix a torus invariant closed locus $Y$ of codimension two in $X$.
Then, by the explicit description of the fan of $X$ (see \cite{CLS} Example 7.3.5.) and the Orbit-Cone correspondence,
one can show that $Y$ is also a projective space bundle over a liner subspase $\mathbb{P}^{s'} \subset \mathbb{P}^s$.
More precisely, $Y = \mathbb{P}_{\mathbb{P}^{s'}}(\mathcal{F})$ 
where $\mathcal{F}$ is a direct sum of $r' +1$ line bundles in $\{ \stsh_{\mathbb{P}^{s'}}(a_i) \}$.
In other words, $Y$ is the intersection of two torus invariant divisors in the linear systems $\lvert p^*\stsh_{\mathbb{P}^s}(1) \rvert$, or $\lvert p^*\stsh_{\mathbb{P}^s}(-a_{\lambda}) \otimes \stsh_p(1) \rvert$ for some $a_{\lambda}$s.
Note that $r' + s' = r + s -2$.

\[ \begin{tikzcd}
    \widetilde{X} \arrow{r}{f}  & X \arrow{r}{p} & \mathbb{P}^s \\
    E \arrow[hook]{u}{\iota} \arrow{r}{\pi} & Y \arrow[hook]{u}  \arrow{r}{q} & \mathbb{P}^{s'} \arrow[hook]{u}  
\end{tikzcd}
\]

\subsection{Mutations} \label{sec:pf (i) mutation}

By Orlov's blow-up formula \ref{blow up}, we obtain a following semiorthogonal decomposition
\[ \D(\widetilde{X}) = \langle \iota_*\pi^*\D(Y) \otimes \stsh(E), f^*\D(X) \rangle. \]
By Theorem \ref{projective bundle}, $\D(X)$ and $\D(Y)$ have exceptional collections
\begin{align*}
\D(X) &= \langle \mathcal{A}, \mathcal{A} \otimes \stsh_p(1), \dots, \mathcal{A} \otimes \stsh_p(r) \rangle, \\
\D(Y) &= \langle\mathcal{A}', \mathcal{A}' \otimes \stsh_q(1), \dots, \mathcal{A}' \otimes \stsh_q(r') \rangle,
\end{align*}
where 
\begin{align*}
\mathcal{A} &= p^*\D(\mathbb{P}^s) = \langle p^*\stsh, p^*\stsh(1), \dots, p^*\stsh(s) \rangle, \\
\mathcal{A}' &= q^*\D(\mathbb{P}^{s'}) = \langle q^*\stsh, q^*\stsh(1), \dots, q^*\stsh(s') \rangle.
\end{align*}
We note that these full exceptional collections in $\D(X)$ and $\D(Y)$ are strong since the bundle $\mathcal{E}$ (resp. $\mathcal{F}$) splits into non-negative line bundles on $\mathbb{P}^s$ (resp. $\mathbb{P}^{s'}$).

In the following, we arrange the pair of integers $(\alpha, \beta)$ in reverse lexicographic order.
This means, we define $(\alpha_1, \beta_1) < (\alpha_2, \beta_2)$ if $\beta_1 < \beta_2$, or $\beta_1 = \beta_2$ and $\alpha_1 < \alpha_2$.

For sake of simplicity, we denote the sheaves on $\widetilde{X}$ by
\begin{align*}
\mathcal{L}_{\alpha, \beta} &:= f^*(p^*\stsh_{\mathbb{P}^s}(\alpha) \otimes \stsh_p(\beta)), \\
\mathcal{M}_{\alpha, \beta} &:= \iota_*\pi^*(q^*\stsh_{\mathbb{P}^{s'}}(\alpha) \otimes \stsh_q(\beta)) \otimes \stsh(E).
\end{align*}

By Lemma \ref{lemA} and Lemma \ref{lem mutation}(i), for the exceptional pair $(\mathcal{M}_{\alpha_1, \beta_1}, \mathcal{L}_{\alpha_2, \beta_2})$ with $(0,0) \leq (\alpha_2, \beta_2) < (\alpha_1, \beta_1) \leq (s', r')$, the right mutation does not change $\mathcal{M}_{\alpha_1, \beta_1}$, i.e.
\[ \mathbb{R}_{\mathcal{L}_{\alpha_2, \beta_2}}(\mathcal{M}_{\alpha_1, \beta_1}) = \mathcal{M}_{\alpha_1, \beta_1}. \]

In addition, if $(\alpha_2, \beta_2) = (\alpha_1, \beta_1)$, we can compute the right mutation as below.

\begin{claim}
For $(\alpha, \beta) \leq (s', r')$, the right mutation for the exceptional pair $(\mathcal{M}_{\alpha, \beta}, \mathcal{L}_{\alpha, \beta})$ is given by
\[ \mathbb{R}_{\mathcal{L}_{\alpha, \beta}}(\mathcal{M}_{\alpha, \beta}) = \mathcal{L}_{\alpha, \beta} \otimes \stsh(E). \]
\end{claim}

From now on, we denote this line bundle by $\mathcal{L}'_{\alpha, \beta} := \mathcal{L}_{\alpha, \beta} \otimes \stsh(E)$. 

\begin{pf} \rm
By Lemma \ref{lemA}, we have
\[ \RHom(\mathcal{M}_{\alpha, \beta}, \mathcal{L}_{\alpha, \beta}) \simeq  \mathbb{C}[-1]. \]
Hence the exact triangle that defines the right mutation
\begin{align*}
\mathbb{R}_{\mathcal{L}_{\alpha, \beta}}(\mathcal{M}_{\alpha, \beta}) \to \mathcal{M}_{\alpha, \beta} \to \mathcal{L}_{\alpha, \beta}[1] 
\end{align*}
coincides with the 1-shifted fundamental sequence
\begin{align*}
\mathcal{L}'_{\alpha, \beta} \to \mathcal{M}_{\alpha, \beta} \to \mathcal{L}_{\alpha, \beta}[1],
\end{align*}
and the uniqueness of mapping cone implies the isomorphism we want. \qed
\end{pf}

Now we apply a mutation operation to above full exceptional collection in order to construct a full exceptional collection consisting of line bundles. First, we right-mutate $\mathcal{M}_{s', r'}$ through objects $\mathcal{L}_{0,0}, \dots, \mathcal{L}_{s'-1, r'}$. 

\begin{center}
\begin{tikzpicture}
 \matrix[matrix of math nodes, column sep=.3em]{
  {{\mathcal{M}}_{0,0}} & \cdots & \mathcal{M}_{s'-1,r'} & |(Msr)[rectangle,draw=black]| \mathcal{M}_{s',r'} & 
  \mathcal{L}_{0,0} & \cdots & |(Lsrm1)| \mathcal{L}_{s'-1,r'} & |(Lsr)| \mathcal{L}_{s',r'} & \mathcal{L}_{s'+1,r'} & \cdots \\};
 \path[->,bend right] (Msr.south) edge ($(Lsrm1.south)!.5!(Lsr.south)$);
\end{tikzpicture}
\end{center}

This mutation does not change $\mathcal{M}_{s', r'}$. Next, we right-mutate $\mathcal{M}_{s', r'}$ through $\mathcal{L}_{s', r'}$.

\begin{center}
\begin{tikzpicture}
\matrix[matrix of math nodes, column sep=.3em]{
\mathcal{M}_{0,0} & \cdots & \mathcal{M}_{s'-1,r'} & \mathcal{L}_{0,0} & \cdots & \mathcal{L}_{s'-1,r'} & |(Msr)[rectangle,draw=black]| \mathcal{M}_{s',r'} & |(Lsrm1)| \mathcal{L}_{s',r'} & |(Lsr)| \mathcal{L}_{s'+1,r'} & \cdots \\};
 \path[->,bend right] (Msr.south) edge ($(Lsrm1.south)!.5!(Lsr.south)$);
\end{tikzpicture}
\end{center}

Then, we have an exceptional collection

\begin{center}
\begin{tikzpicture}
\matrix[matrix of math nodes, column sep=.2em]{
\mathcal{M}_{0,0} & \cdots & \mathcal{M}_{s'-1,r'} & \mathcal{L}_{0,0} & \cdots & \mathcal{L}_{s'-1,r'} & \mathcal{L}_{s',r'} & |[rectangle,draw=black]| \mathcal{L}'_{s',r'} & \mathcal{L}_{s'+1,r'} & \cdots \\};
\end{tikzpicture}
\end{center}

In the same way as above, we apply the mutation operations for $\mathcal{M}_{s'-1, r'}$, $\mathcal{M}_{s'-2,r'}$, $\dots$, $\mathcal{M}_{0,0}$ one after the other. After this operation, we finally obtain the full exceptional collections consisting of line bundles $\{\mathcal{L}_{\alpha, \beta}\}_{\alpha, \beta}$ with $0 \leq \alpha \leq s, ~0 \leq \beta \leq r$ and $\{\mathcal{L}'_{\alpha, \beta}\}_{\alpha, \beta}$ with $0 \leq \alpha \leq s', ~0 \leq \beta \leq r'$ orderd by $ \mathcal{L}_{\alpha_1, \beta_1} \leq \mathcal{L}'_{\alpha_1, \beta_1} \leq \mathcal{L}_{\alpha_2, \beta_2}$ for $(\alpha_1, \beta_1) < (\alpha_2, \beta_2)$.

\[ \D(\widetilde{X}) = \langle \mathcal{L}_{0,0}, ~\mathcal{L}'_{0,0}, ~\mathcal{L}_{1,0}, ~\mathcal{L}'_{1,0}, ~\cdots, ~\mathcal{L}_{s',r'} , ~\mathcal{L}'_{s',r'}, ~\mathcal{L}_{s'+1,r'}, ~\cdots, ~\mathcal{L}_{s,r} \rangle. \]

\subsection{Strongness}

In this subsection, we write $\stsh(\alpha)$ instead of $p^*\stsh_{\mathbb{P}^s}(\alpha)$. 
The aim of this subsection is to prove the following lemma.

\begin{lem} \label{lem 5.2}
The exceptional collection of line bundles which we constructed in the above subsection is strong.
\end{lem}

\begin{pf}\rm
What is nontrivial is the following vanishing and other vanishing of extensions we need follows from Lemma \ref{lemB}.
\begin{align*}
&\Ext^i_{\widetilde{X}}(\mathcal{L}'_{\alpha_1, \beta_1}, \mathcal{L}_{\alpha_2, \beta_2}) \\
= &\Ext^i_{\widetilde{X}}(f^*(\stsh(\alpha_1) \otimes \stsh_p(\beta_1)) \otimes \stsh(E), f^*(\stsh(\alpha_2) \otimes \stsh_p(\beta_2))) \\
= &0 \\
&\text{for all $i \neq 0$ and $\left\{
\begin{array}{ll}
0 \leq \beta_1 < \beta_2 \leq r ~(\beta_1 \leq r'), 0 \leq \alpha_1 \leq s', 0 \leq \alpha_2 \leq s, \\
\text{or}~~ \beta_1 = \beta_2, 0 \leq \alpha_1 < \alpha_2 \leq s ~( \alpha_1 \leq s').
\end{array}
\right.$} 
\end{align*}
By using the projection formula, we have an isomorphism
\begin{align*}
&\Ext^i_{\widetilde{X}}(f^*(\stsh(\alpha_1) \otimes \stsh_p(\beta_1)) \otimes \stsh(E), f^*(\stsh(\alpha_2) \otimes \stsh_p(\beta_2))) \\
\simeq ~ &\Ho^i(X, \stsh(\alpha_2 - \alpha_1) \otimes \stsh_p(\beta_2 - \beta_1) \otimes I_Y)
\end{align*}
A short exact sequence on $X$
\begin{align*}
0 \to \stsh(\alpha_2 - \alpha_1) \otimes \stsh_p(\beta_2 - \beta_1) &\otimes I_Y \to \stsh(\alpha_2 - \alpha_1) \otimes \stsh_p(\beta_2 - \beta_1) \\
&\to (\stsh(\alpha_2 - \alpha_1) \otimes \stsh_p(\beta_2 - \beta_1))|_{Y} \to 0
\end{align*}
and the vanishing of cohomologies 
\begin{align*}
\Ho^i(X, \stsh(\alpha_2 - \alpha_1) \otimes \stsh_p(\beta_2 - \beta_1)) &= 0, \\
\Ho^i(Y, \stsh(\alpha_2 - \alpha_1) \otimes \stsh_q(\beta_2 - \beta_1)) &= 0
\end{align*}
for all $i > 0$ imply that
\[ \Ho^i(X, \stsh(\alpha_2 - \alpha_1) \otimes \stsh_p(\beta_2 - \beta_1) \otimes I_Y) = 0 \]
for all $i \geq 2$. 
To prove the vanishing
\[ \Ho^1(X, \stsh(\alpha_2 - \alpha_1) \otimes \stsh_p(\beta_2 - \beta_1) \otimes I_Y) = 0, \]
we need to check the surjectivity of the map
\[ \Ho^0(X, \stsh(\alpha_2 - \alpha_1) \otimes \stsh_p(\beta_2 - \beta_1)) \to \Ho^0(Y, \stsh(\alpha_2 - \alpha_1) \otimes \stsh_q(\beta_2 - \beta_1)). \]
This is equivalent to the surjecvity of
\[ \Ho^0(\mathbb{P}^s, \stsh(\alpha_2 - \alpha_1) \otimes \mathrm{Sym}^{\beta_2 - \beta_1}\mathcal{E}) \to \Ho^0(\mathbb{P}^{s'}, \stsh(\alpha_2 - \alpha_1) \otimes \mathrm{Sym}^{\beta_2 - \beta_1}\mathcal{F}). \]
Because the bundle $\mathcal{E}$ splits into a direct sum of positive line bundles on $\mathbb{P}^s$ and $\alpha_2 - \alpha_1 \geq -s'$, the restriction morphism
\[ \Ho^0(\mathbb{P}^s, \stsh(\alpha_2 - \alpha_1) \otimes \mathrm{Sym}^{\beta_2 - \beta_1}\mathcal{E}) \twoheadrightarrow \Ho^0(\mathbb{P}^{s'}, \stsh(\alpha_2 - \alpha_1) \otimes \mathrm{Sym}^{\beta_2 - \beta_1}\mathcal{E}'). \]
is surjective. Furthermore, since $\mathcal{E}' = \mathcal{E}|_{\mathbb{P}^{s'}}$ splits as $\mathcal{E}' = \mathcal{G} \oplus \mathcal{F}$, the morphism
\[ \Ho^0(\mathbb{P}^{s'}, \stsh(\alpha_2 - \alpha_1) \otimes \mathrm{Sym}^{\beta_2 - \beta_1}\mathcal{E}') \to \Ho^0(\mathbb{P}^{s'}, \stsh(\alpha_2 - \alpha_1) \otimes \mathrm{Sym}^{\beta_2 - \beta_1}\mathcal{F}). \]
coincieds the projection morphism, and hence is also surjective. Thus, the proof was completed.
\qed
\end{pf}

\section{Proof of Theorem \ref{Main Theorem}, codimension three case}

Let $\mathcal{E} = \stsh_{\mathbb{P}^s}(a_0) \oplus \stsh_{\mathbb{P}^s}(a_1) \oplus \cdots \oplus \stsh_{\mathbb{P}^s}(a_r) ~ (0 = a_0 \leq a_1 \leq \cdots \leq a_r)$ a vector bundles on $\mathbb{P}^s$ such that $X = \mathbb{P}_{\mathbb{P}^s}(\mathcal{E})$, 
and we set $a = \sum_{k=0}^r a_k$.
As in the above section, we can set $Y = \mathbb{P}_{\mathbb{P}^{s'}}(\mathcal{F})$ where $\mathcal{F}$ is a direct sum of $r' +1$ line bundles in $\{ \stsh_{\mathbb{P}^{s'}}(a_i) \}_{i=0}^r$. Note that $r' + s' = r + s - 3$.
In other words, $Y$ is the intersection of three torus invariant divisors in the linear systems $\lvert p^*\stsh_{\mathbb{P}^s}(1) \rvert$,
or $\lvert p^*\stsh_{\mathbb{P}^s}(-a_{\lambda}) \otimes \stsh_p(1) \rvert$ for some $a_{\lambda}$s.

Note that the canonical bundle of $\widetilde{X}$ is given by
\begin{align*}
 \omega_{\widetilde{X}} &= f^*\omega_X \otimes \stsh(2E) \\
&= f^*p^*\stsh_{\mathbb{P}^s}(-s-1+a) \otimes f^*\stsh_p(-r-1) \otimes \stsh(2E).
\end{align*}

\[ \begin{tikzcd}
    \widetilde{X} \arrow{r}{f}  & X \arrow{r}{p} & \mathbb{P}^s \\
    E \arrow[hook]{u}{\iota} \arrow{r}{\pi} & Y \arrow[hook]{u}  \arrow{r}{q} & \mathbb{P}^{s'} \arrow[hook]{u}  
\end{tikzcd}
\]

\subsection{Mutations}

By Orlov's blow-up formula \ref{blow up}, we obtain the following semiorthogonal decomposition
\[ \D(\widetilde{X}) = \langle \iota_*\pi^*\mathcal{D}_2 \otimes \stsh(2E), \iota_*\pi^*\mathcal{D}_1 \otimes \stsh(E), f^*\D(X) \rangle, \]
where $\mathcal{D}_1 = \mathcal{D}_2 = \D(Y)$.
By Lemma \ref{lem mutation}, we mutate $\iota_*\pi^*\mathcal{D}_2 \otimes \stsh(2E)$ and obtain another semiorthogonal decomposition of $\D(\widetilde{X})$,
\begin{align*}
\D(\widetilde{X}) &= \langle \iota_*\pi^*\mathcal{D}_1 \otimes \stsh(E), f^*\D(X), \iota_*\pi^*\mathcal{D}_2 \otimes f^*\omega_X^{-1} \rangle, \\
&= \langle \iota_*\pi^*\mathcal{D}_1 \otimes \stsh(E), f^*\D(X), \iota_*\pi^*(\mathcal{D}_2 \otimes \omega_X^{-1}|_{Y}) \rangle.
\end{align*}
The derived category $\D(X)$ of $X$ has an exceptional collection
\[ \D(X) = \langle \mathcal{A}, \mathcal{A} \otimes \stsh_p(1), \dots, \mathcal{A} \otimes \stsh_p(r) \rangle, \]
where 
\[ \mathcal{A} = p^*\D(\mathbb{P}^s) = \langle p^*\stsh, p^*\stsh(1), \dots, p^*\stsh(s) \rangle. \]
We take full exceptional collections of the categories $\mathcal{D}_1$ and $\mathcal{D}_2$ as
\begin{align*}
\mathcal{D}_1 &= \langle\mathcal{A}', \mathcal{A}' \otimes \stsh_q(1), \dots, \mathcal{A}' \otimes \stsh_q(r') \rangle, \\
\mathcal{D}_2 &= \langle\mathcal{A}'' \otimes \stsh_q(-r'-1), \mathcal{A}'' \otimes \stsh_q(-r'), \dots, \mathcal{A}'' \otimes \stsh_q(-1) \rangle,
\end{align*}
where
\begin{align*}
\mathcal{A}' &= q^*\D(\mathbb{P}^{s'}) = \langle q^*\stsh, q^*\stsh(1), \dots, q^*\stsh(s') \rangle, ~~\text{and}\\
\mathcal{A}'' &= q^*\D(\mathbb{P}^{s'}) = \langle q^*\stsh(-s'-1+a), q^*\stsh(-s'+a), \dots, q^*\stsh(-1+a) \rangle,
\end{align*}
respectively. Then, we have
\begin{align*}
\mathcal{D}_2 \otimes \omega_X^{-1}|_Y &= \langle \mathcal{A}''  \otimes q^*\stsh(s+1-a)  \otimes \stsh_q(r-r'), \\
&\mathcal{A}'' \otimes q^*\stsh(s+1-a) \otimes \stsh_q(r-r'+1), \dots, \mathcal{A}'' \otimes q^*\stsh(s+1-a) \otimes \stsh_q(r) \rangle,
\end{align*}
and
\begin{align*}
\mathcal{A}'' \otimes q^*\stsh(s+1-a) = q^*\D(\mathbb{P}^{s'}) = \langle q^*\stsh(s-s'), q^*\stsh(s-s'+1), \dots, q^*\stsh(s) \rangle.
\end{align*}

We apply exactly the same sequence of mutations as in Section \ref{sec:pf (i) mutation} to the part $\langle \iota_*\pi^*\mathcal{D}_1 \otimes \stsh(E), \pi^*\D(X) \rangle$, and obtain the following exceptional collection
\begin{align*}
\D(\widetilde{X}) = \langle \mathcal{B},  & \mathcal{B} \otimes f^*\stsh_p(1), \cdots, \mathcal{B} \otimes f^*\stsh_p(r'), \\
&f^*\mathcal{A} \otimes f^*\stsh_p(r'+1), \cdots, f^*\mathcal{A} \otimes f^*\stsh_p(r), \iota_*\pi^*(\mathcal{D}_2 \otimes \omega_X^{-1}|_Y) \rangle,
\end{align*}
where
\begin{align*}
\mathcal{B} = \langle \stsh, \stsh(E), &f^*p^*\stsh(1), f^*p^*\stsh(1) \otimes \stsh(E), \cdots, \\
& f^*p^*\stsh(s'), f^*p^*\stsh(s') \otimes \stsh(E), f^*p^*\stsh(s'+1), \cdots, f^*p^*\stsh(s) \rangle.
\end{align*}

In the following, we denote the sheaves on $\widetilde{X}$ by
\begin{align*}
\mathcal{L}_{\alpha, \beta} &:= f^*(p^*\stsh_{\mathbb{P}^s}(\alpha) \otimes \stsh_p(\beta)), \\
\mathcal{L}'_{\alpha, \beta} &:= \mathcal{L}_{\alpha, \beta} \otimes \stsh(E), \\
\mathcal{L}''_{\alpha, \beta} &:= \mathcal{L}_{\alpha, \beta} \otimes \stsh(-E), \\
\mathcal{M}'_{\alpha, \beta} &:= \iota_*\pi^*(q^*\stsh_{\mathbb{P}^{s'}}(\alpha) \otimes \stsh_q(\beta))
\end{align*}
for brevity.

Next, we mutate the exceptional objects in $\iota_*\pi^*(\mathcal{D}_2 \otimes \omega_X^{-1}|_Y)$. 
In order to compute the mutations explicitly, we need the following lemma.

\begin{lem} \label{lem 5-1}
The following holds.
\begin{enumerate}
\item[(a)] The extensions of sheaves on $\widetilde{X}$
\[ \Ext^i_{\widetilde{X}}(\mathcal{L}_{\alpha_1, \beta_1}, \mathcal{M}_{\alpha_2, \beta_2})  \]
is zero for all $i \in \mathbb{Z}$ and 
\[ \left\{
\begin{array}{ll}
0 \leq \beta_1 \leq r, r-r' \leq \beta_2 \leq r, \beta_2 < \beta_1, 0 \leq \alpha_1 \leq s, s-s' \leq \alpha_2 \leq s,\\
\text{or}~~ \beta_1 = \beta_2, 0 \leq \alpha_2 < \alpha_1 \leq s.
\end{array}
\right. \]

\item[(b)] The extensions of sheaves on $\widetilde{X}$
\[ \Ext^i_{\widetilde{X}}(\mathcal{L}'_{\alpha_1, \beta_1}, \mathcal{M}_{\alpha_2, \beta_2}) \] 
is zero for all $i \in \mathbb{Z}$ and 
\[ \left\{
\begin{array}{ll}
0 \leq \beta_1 \leq r', r-r' \leq \beta_2 \leq r, \beta_2 < \beta_1, 0 \leq \alpha_1 \leq s', s-s' \leq \alpha_2 \leq s, \\
\text{or}~~ \beta_1 = \beta_2, 0 \leq \alpha_1 \leq s', s-s' \leq \alpha_2 \leq s, \alpha_2 \leq \alpha_1.
\end{array}
\right. \]
\end{enumerate}
\end{lem}

\begin{pf}
\rm
For (a), we have
\begin{align*}
&\Ext^i_{\widetilde{X}}(\mathcal{L}_{\alpha_1, \beta_1}, \mathcal{M}_{\alpha_2, \beta_2}) \\
= &\Ext^i_{\widetilde{X}}(f^*p^*\stsh_{\mathbb{P}^s}(\alpha_1) \otimes f^*\stsh_p(\beta_1), \iota_*\pi^*(q^*\stsh_{\mathbb{P}^{s'}}(\alpha_2) \otimes \stsh_q(\beta_2))) \\
\simeq ~ &\Ho^i(Y, q^*\stsh_{\mathbb{P}^{s'}}(\alpha_2 - \alpha_1) \otimes \stsh_q(\beta_2 - \beta_1)) \\
\simeq ~ &\Ho^i(\mathbb{P}^{s'}, \stsh_{\mathbb{P}^{s'}}(\alpha_2 - \alpha_1) \otimes Rq_*\stsh_q(\beta_2 - \beta_1)) \\
= ~ &0
\end{align*}
for $(\alpha_1, \beta_1)$ and $(\alpha_2, \beta_2)$ that satisfy the above condition.

For (b), first, we have
\begin{align*}
&\Ext^i_{\widetilde{X}}(\mathcal{L}'_{\alpha_1, \beta_1}, \mathcal{M}_{\alpha_2, \beta_2}) \\
= &\Ext^i_{\widetilde{X}}(f^*p^*\stsh_{\mathbb{P}^s}(\alpha_1) \otimes f^*\stsh_p(\beta_1) \otimes \stsh(E), \iota_*\pi^*(q^*\stsh_{\mathbb{P}^{s'}}(\alpha_2) \otimes \stsh_q(\beta_2))) \\
\simeq ~ &\Ho^i(E, \pi^*(q^*\stsh_{\mathbb{P}^{s'}}(\alpha_2 - \alpha_1) \otimes \stsh_q(\beta_2 - \beta_1)) \otimes \stsh_E(-E)) \\
\simeq ~ &\Ho^i(Y, q^*\stsh_{\mathbb{P}^{s'}}(\alpha_2 - \alpha_1) \otimes \stsh_q(\beta_2 - \beta_1) \otimes \mathcal{N}_{Y/X}^*).
\end{align*}
The conormal bundle $\mathcal{N}_{Y/X}^*$ of $Y \subset X$ splits into three line bundles each of which is of the form $q^*\stsh_{\mathbb{P}^{s'}}(-1)$ or $q^*\stsh_{\mathbb{P}^{s'}}(a_{\lambda}) \otimes \stsh_q(-1)$ for some $0 \leq \lambda \leq r$.
From now on, we check the vanishing of this cohomology.
Here we prove this only in the case $r' = r - 2$ and $s' = s-1$, 
but the reader can easily prove other cases by the same argument.

In this case, the conormal bundle of $Y$ is given by
\[ \mathcal{N}_{Y/X}^* = q^*\stsh_{\mathbb{P}^{s'}}(-1) \oplus \bigoplus_{k=1}^2 (q^*\stsh_{\mathbb{P}^{s'}}(a_{\lambda_k}) \otimes \stsh_q(-1)). \]
Then, we have
\begin{align*}
\Ho^i(Y, q^*\stsh_{\mathbb{P}^{s'}}(\alpha_2 - \alpha_1 -1) \otimes \stsh_q(\beta_2 - \beta_1)) = 0
\end{align*}
for all $i \in \mathbb{Z}$, since $- r' \leq \beta_2 - \beta_1 <0$, or $\beta_2=\beta_1$ and $0 > \alpha_2 - \alpha_1 -1 \geq s-s'-s'-1 = -s'$, and we have
\begin{align*}
\Ho^i(Y, q^*\stsh_{\mathbb{P}^{s'}}(\alpha_2 - \alpha_1 + a_{\lambda_k}) \otimes \stsh_q(\beta_2 - \beta_1 - 1)) = 0
\end{align*}
for all $i \in \mathbb{Z}$, since $0 > \beta_2 - \beta_1 - 1 \geq r - r' - r' -1 = -r' + 1 > - r'$. 
Hence we have the desired vanishing of cohomologies. \qed
\end{pf}

First, we left-mutate $\mathcal{M}'_{s-s', r-r'}$. By Lemma \ref{lem 5-1} and Lemma \ref{lem mutation}, the left mutations of $\mathcal{M}'_{s-s', r-r'}$ over line bundles $\mathcal{L}_{s,r}$, $\dots$, $\mathcal{L}'_{s-s', r-r'}$ do not change $\mathcal{M}'_{s-s',r-r'}$.

\begin{center}
\begin{tikzpicture}
 \matrix[matrix of math nodes, column sep=.3em]{
 \cdots & \mathcal{L}'_{s-s'-1,r-r'} & |(Lsrm1)| \mathcal{L}_{s-s',r-r'} & |(Lsr)| \mathcal{L}'_{s-s',r-r'} & \cdots & \mathcal{L}_{s,r} & |(Msr)[rectangle,draw=black]| \mathcal{M}'_{s-s',r-r'} & \cdots &  \mathcal{M}'_{s,r} \\};
 \path[->,bend left] (Msr.south) edge ($(Lsrm1.south)!.5!(Lsr.south)$);
\end{tikzpicture}
\end{center}

Next, we left-mutate $\mathcal{M}_{s-s',r-r'}$ over $\mathcal{L}_{s-s', r-r'}$.

\begin{center}
\begin{tikzpicture}
 \matrix[matrix of math nodes, column sep=.3em]{
 \cdots & |(Lsrm1)| \mathcal{L}'_{s-s'-1,r-r'} & |(Lsr)| \mathcal{L}_{s-s',r-r'} & |(Msr)[rectangle,draw=black]| \mathcal{M}'_{s-s',r-r'} & \mathcal{L}'_{s-s',r-r'} & \cdots & \mathcal{L}_{s,r}  & \cdots &  \mathcal{M}'_{s,r} \\};
 \path[->,bend left] (Msr.south) edge ($(Lsrm1.south)!.5!(Lsr.south)$);
\end{tikzpicture}
\end{center}

By the isomorphism $\mathbb{L}_{\mathcal{L}_{\alpha, \beta}}(\mathcal{M}'_{\alpha, \beta}) \simeq \mathcal{L}_{\alpha, \beta} \otimes \stsh(-E) =: \mathcal{L}''_{\alpha, \beta}$, we have a new exceptional collection

\begin{center}
\begin{tikzpicture}
 \matrix[matrix of math nodes, column sep=.3em]{
 \cdots &\mathcal{L}'_{s-s'-1,r-r'} & \mathcal{L}''_{s-s',r-r'} & \mathcal{L}_{s-s',r-r'} & \mathcal{L}'_{s-s',r-r'} & \cdots & \mathcal{L}_{s,r}  & \cdots &  \mathcal{M}'_{s,r}. \\};
\end{tikzpicture}
\end{center}

In the same way as above, we apply mutation operations to, $\mathcal{M}_{s-s'+1, r-r'}$, $\mathcal{M}_{s-s'+2, r-r'}$, $\dots$, $\mathcal{M}_{s,r}$ one after another. Finally, we get a full exceptional collection consisting of line bundles
\begin{align*}
&\{\mathcal{L}_{\alpha, \beta} \}_{\alpha, \beta}  ~\text{with $0 \leq \alpha \leq s, ~0 \leq \beta \leq r$} \\
&\{\mathcal{L}'_{\alpha, \beta} \}_{\alpha, \beta}  ~\text{with $0 \leq \alpha \leq s', ~0 \leq \beta \leq r'$} \\
&\{\mathcal{L}''_{\alpha, \beta} \}_{\alpha, \beta}  ~\text{with $s-s' \leq \alpha \leq s, ~r-r' \leq \beta \leq r$}.
\end{align*}
placed in ascending order defined by $\mathcal{L}_{\alpha_1, \beta_1} \leq \mathcal{L}''_{\alpha_2, \beta_2} \leq \mathcal{L}_{\alpha_2, \beta_2} \leq \mathcal{L}'_{\alpha_2, \beta_2} \leq \mathcal{L}_{\alpha_3, \beta_3}$ for $(\alpha_1, \beta_1) < (\alpha_2, \beta_2) < (\alpha_3, \beta_3)$.

\subsection{Strongness}

\begin{lem}
The full exceptional collection of line bundles which is constructed in the above subsection is strong.
\end{lem}

\begin{pf} \rm
What is non-trivial is to show that the following vanishings and other vanishings we need follow from Lemma \ref{lemB}.
\begin{enumerate}
\item[(A)] $\Ext^i_{\widetilde{X}}(\mathcal{L}_{\alpha_1, \beta_1}, \mathcal{L}''_{\alpha_2, \beta_2}) = 0$
\[ \text{for all $i \neq 0$ and $\left\{
\begin{array}{ll}
0 \leq \beta_1 < \beta_2 \leq r', 0 \leq \alpha_1 \leq s, s-s' \leq \alpha_2 \leq s, \\
\text{or}~~ \beta_1 = \beta_2, 0 \leq \alpha_1 \leq s, s-s' \leq \alpha_2 \leq s, \alpha_1 < \alpha_2.
\end{array}
\right.$} \]

\item[(B)] $\Ext^i_{\widetilde{X}}(\mathcal{L}'_{\alpha_1, \beta_1}, \mathcal{L}''_{\alpha_2, \beta_2}) = 0$
\[ \text{for all $i \neq 0$ and $\left\{
\begin{array}{ll}
0 \leq \beta_1 \leq r', r-r' \leq \beta_2 \leq r, \beta_1 < \beta_2, 0 \leq \alpha_1 \leq s', s-s' \leq \alpha_2 \leq s, \\
\text{or}~~ \beta_1 = \beta_2, 0 \leq \alpha_1 \leq s', s-s' \leq \alpha_2 \leq s, \alpha_1 < \alpha_2.
\end{array}
\right.$} \]

\item[(C)] $\Ext^i_{\widetilde{X}}(\mathcal{L}'_{\alpha_1, \beta_1}, \mathcal{L}_{\alpha_2, \beta_2}) = 0$
\[ \text{for all $i \neq 0$ and $\left\{
\begin{array}{ll}
0 \leq \beta_1 \leq r', 0 \leq \beta_2 \leq r, \beta_1 < \beta_2, 0 \leq \alpha_1 \leq s', 0 \leq \alpha_2 \leq s, \\
\text{or}~~ \beta_1 = \beta_2, 0 \leq \alpha_1 \leq s', 0 \leq \alpha_1 < \alpha_2 \leq s.
\end{array}
\right.$} \]
\end{enumerate}
For the first, we have an isomorphism
\begin{align*}
&\Ho^i(\widetilde{X}, f^*(p^*\stsh_{\mathbb{P}^s}(\alpha) \otimes \stsh_p(\beta)) \otimes \stsh(-E)) \\
\simeq ~&\Ho^i(X, p^*\stsh_{\mathbb{P}^s}(\alpha) \otimes \stsh_p(\beta) \otimes I_Y).
\end{align*}
Let us consider the exact sequence
\[ 0 \to p^*\stsh_{\mathbb{P}^s}(\alpha) \otimes \stsh_p(\beta) \otimes I_Y \to p^*\stsh_{\mathbb{P}^s}(\alpha) \otimes \stsh_p(\beta) \to (p^*\stsh_{\mathbb{P}^s}(\alpha) \otimes \stsh_p(\beta))|_{Y} \to 0. \]
The cohomologies of the second and third terms vanish
\begin{align*}
&\Ho^i(X, p^*\stsh_{\mathbb{P}^s}(\alpha) \otimes \stsh_p(\beta)) = 0, \\
&\Ho^i(Y, q^*\stsh_{\mathbb{P}^{s'}}(\alpha) \otimes \stsh_q(\beta)) = 0
\end{align*}
for all $i > 0$ and for all $\beta \geq 0$ and $\alpha \geq -s'$. By combining it with the subjectivity of the map
\[ \Ho^i(X, p^*\stsh_{\mathbb{P}^s}(\alpha) \otimes \stsh_p(\beta)) \twoheadrightarrow \Ho^i(Y, q^*\stsh_{\mathbb{P}^{s'}}(\alpha) \otimes \stsh_q(\beta)) \]
that follows from the same argument as in the last part of the proof of Lemma \ref{lem 5.2},
we have
\[ \Ho^i(X, p^*\stsh_{\mathbb{P}^s}(\alpha) \otimes \stsh_p(\beta) \otimes I_Y) = 0 \]
for all $i \neq 0$ and for all $\beta \geq 0$ and $\alpha \geq -s'$. This proves (A) and (C).

It remains to show (B). 
First, we have an isomorphism
\begin{align*}
&\Ho^i(\widetilde{X}, f^*(p^*\stsh_{\mathbb{P}^s}(\alpha) \otimes \stsh_p(\beta)) \otimes \stsh(-2E)) \\
\simeq ~&\Ho^i(X, p^*\stsh_{\mathbb{P}^s}(\alpha) \otimes \stsh_p(\beta) \otimes I_Y^2).
\end{align*}
Let us consider the exact sequence
\[ 0 \to p^*\stsh_{\mathbb{P}^s}(\alpha) \otimes \stsh_p(\beta) \otimes I_Y^2 \to p^*\stsh_{\mathbb{P}^s}(\alpha) \otimes \stsh_p(\beta) \otimes I_Y \to q^*\stsh_{\mathbb{P}^{s'}}(\alpha) \otimes \stsh_q(\beta) \otimes \mathcal{N}_{Y/X}^* \to 0. \]
It follows from the above computation that the cohomology of the second term vanishes:
\[ \Ho^i(X, p^*\stsh_{\mathbb{P}^s}(\alpha) \otimes \stsh_p(\beta) \otimes I_Y) = 0 \]
for all $i \neq 0$ and for all $\beta \geq 0$ and $\alpha \geq s-2s' \geq -s'$. 

Next, we calculate the cohomology of the third term.
As we proved in the proof of Lemma \ref{lem 5-1} (b),
\[ \Ho^i(Y, q^*\stsh_{\mathbb{P}^{s'}}(\alpha) \otimes \stsh_q(\beta) \otimes \mathcal{N}_{Y/X}^*) = 0 \]
for all $i \neq 0$, and consequently we have
\[ \Ho^i(X, p^*\stsh_{\mathbb{P}^s}(\alpha) \otimes \stsh_p(\beta) \otimes I_Y^2) = 0 \]
for all $i \neq 0, 1$. In order to prove the vanishing of $\Ho^1(X, p^*\stsh_{\mathbb{P}^s}(\alpha) \otimes \stsh_p(\beta) \otimes I_Y^2)$, we have to show that the map
\[ \Ho^0(X, p^*\stsh_{\mathbb{P}^s}(\alpha) \otimes \stsh_p(\beta) \otimes I_Y) \to \Ho^0(Y, q^*\stsh_{\mathbb{P}^{s'}}(\alpha) \otimes \stsh_q(\beta) \otimes \mathcal{N}_{Y/X}^*). \]
is surjective.
Let us take torus invariant prime divisors $D_1, D_2, D_3$ on $X$ such that $Y = D_1 \cap D_2 \cap D_3$, and let $Y_1 = D_1$, $Y_2 = D_1 \cap D_2$, and $Y_3 =Y$. In the below, we treat the case $s' = s-1$ and $r' = r-2$. In this case, we can take divisors as $D_i \in \lvert p^*\stsh_{\mathbb{P}^s}(a_{\lambda_i}) \otimes \stsh_p(-1) \rvert ~~(i=1,2)$ and $D_3 \in \lvert p^*\stsh_{\mathbb{P}^s}(-1) \rvert$. In the following, we set 
\[ \mathcal{L} := p^*\stsh_{\mathbb{P}^s}(\alpha) \otimes \stsh_p(\beta). \]

\begin{claim} \label{claim(II)-1} \rm
The map
\[  \Ho^0(X, \mathcal{L} \otimes I_{Y_1}) \to \Ho^0(Y_1, \mathcal{L}|_{Y_1} \otimes \mathcal{N}_{Y_1/X}^*) \]
is surjective, and the first cohomology group of $\mathcal{L} \otimes I_{Y_1}$ isf
\[ \Ho^1(X, \mathcal{L} \otimes I_{Y_1}) = 0. \]
\end{claim}

\begin{pf} 
\rm
As we have $I_{Y_1} = \stsh_X(-D_1) = p^*\stsh_{\mathbb{P}^s}(a_{\lambda}) \otimes \stsh_p(-1)$, we get
\begin{align*}
\Ho^1(X, \mathcal{L} \otimes I_{Y_1}) &\simeq \Ho^1(X, p^*\stsh_{\mathbb{P}^s}(\alpha + a_{\lambda}) \otimes \stsh_p(\beta-1)) =0,
\end{align*}
and
\begin{align*}
\Ho^1(X, \mathcal{L} \otimes I_{Y_1}^2) &\simeq \Ho^1(X, p^*\stsh_{\mathbb{P}^s}(\alpha + 2a_{\lambda}) \otimes \stsh_p(\beta-2)) = 0,
\end{align*}
since $\beta - 2 \geq -2$ and $\alpha + 2a_{\lambda} \geq -s' = -s$. \qed
\end{pf}

Let us consider a commutative diagram with exact rows:
\[ \begin{tikzcd}
0 \arrow{r} & \mathcal{L} \otimes I_{Y_1} \arrow{r} \arrow{d} & \mathcal{L} \otimes I_{Y_2} \arrow{d} \arrow{r} & \mathcal{L}|_{Y_1} \otimes I_{Y_2/Y_1} \arrow{r} \arrow{d} & 0 \\
0 \arrow{r} & (\mathcal{L} \otimes {\mathcal{N}_{Y_1/X}^*})|_{Y_2} \arrow{r} & \mathcal{L}|_{Y_2} \otimes \mathcal{N}_{Y_2/X}^* \arrow{r} & \mathcal{L}|_{Y_2} \otimes \mathcal{N}_{Y_2/Y_1}^* \arrow{r} & 0
 \end{tikzcd} \]

\begin{claim} \label{claim(II)-2} \rm
The map
\[  \Ho^0(X, \mathcal{L} \otimes I_{Y_2}) \to \Ho^0(Y_2, \mathcal{L}|_{Y_2} \otimes \mathcal{N}_{Y_2/X}^*) \]
is surjective, and we have
\[ \Ho^1(X, \mathcal{L} \otimes I_{Y_2}) = 0. \]
\end{claim}

\begin{pf} \rm
First, we note that 
\begin{align*}
I_{Y_1} &= \stsh_X(-D_1) = p^*\stsh_{\mathbb{P}^s}(a_{\lambda_1}) \otimes \stsh_p(-1), \\
I_{Y_2/Y_1} &= (p^*\stsh_{\mathbb{P}^s}(a_{\lambda_2}) \otimes \stsh_p(-1))|_{Y_1}
\end{align*}
for some $a_{\lambda_1}, a_{\lambda_2} \geq 0$. 
Let us consider the exact sequence
\[ 0 \to \mathcal{L} \otimes I_{Y_1} \to \mathcal{L} \otimes I_{Y_2} \to \mathcal{L}|_{Y_1} \otimes I_{Y_2/Y_1} \to 0. \]
By using the above description of $I_{Y_2/Y_1}$, we have
\begin{align*}
\Ho^1(Y_1, \mathcal{L}|_{Y_1} \otimes I_{Y_2/Y_1}) &\simeq \Ho^1(Y_1, (p^*\stsh_{\mathbb{P}^s}(\alpha + a_{\lambda_2}) \otimes \stsh_p(\beta-1))|_{Y_1}) = 0,
\end{align*}
since $Y_1$ is a $\mathbb{P}^{r-1}$-bundle over $\mathbb{P}^s$ and $\beta - 1 \geq -1$ and $\alpha + a_{\lambda_2} \geq s-2s'+a_{\lambda_2} \geq -s$.
Moreover, by Claim \ref{claim(II)-1}, we have $\Ho^1(X, \mathcal{L} \otimes I_{Y_1}) = 0$,
and hence the vanishing of $\Ho^1(X, \mathcal{L} \otimes I_{Y_2})$ follows.

Next, we prove the surjectivity of the map $\Ho^0(X, \mathcal{L} \otimes I_{Y_2}) \to \Ho^0(Y_2, \mathcal{L}|_{Y_2} \otimes \mathcal{N}_{Y_2/X}^*)$.
First, by Claim \ref{claim(II)-1}, the map
\[ \Ho^0(X, \mathcal{L} \otimes I_{Y_2}) \to \Ho^0(Y_1, \mathcal{L}|_{Y_1} \otimes {I_{Y_2}}_{|Y_1}) \]
is surjective.
We also have
\begin{align*}
\Ho^1(Y_1, \mathcal{L}|_{Y_1} \otimes I_{Y_2/Y_1}^2) &\simeq \Ho^1(Y_1, (p^*\stsh_{\mathbb{P}^s}(\alpha + 2a_{\lambda_2}) \otimes \stsh_p(\beta-2))|_{Y_1}) = 0
\end{align*}
(we note that if $r = 2$ in the Case (B), then $\beta = 2$ by our construction of the full exceptional collection),
and hence the map 
\[ \Ho^0(Y_1, \mathcal{L}|_{Y_1} \otimes I_{Y_2/Y_1}) \to \Ho^0(Y_2, \mathcal{L}|_{Y_2} \otimes \mathcal{N}_{Y_2/Y_1}^*) \]
is surjective.
Next, we consider the exact sequence
\[ 0 \to \mathcal{L}|_{Y_1} \otimes \mathcal{N}_{Y_1/X}^* \otimes I_{Y_2/Y_1} \to \mathcal{L}|_{Y_1} \otimes \mathcal{N}_{Y_1/X}^* \to (\mathcal{L} \otimes {\mathcal{N}_{Y_1/X}^*})|_{Y_2} \to 0. \]
The first cohomology group of the first term of this sequence is
\begin{align*}
\Ho^1(Y_1, \mathcal{L}|_{Y_1} \otimes \mathcal{N}_{Y_1/X}^* \otimes I_{Y_2/Y_1}) &\simeq \Ho^1(Y_1, (\stsh_{\mathbb{P}^s}(\alpha + a_{\lambda_1} + a_{\lambda_2})  \otimes \stsh_p(\beta -2))|_{Y_1}) = 0.
\end{align*}
By combining it with Claim \ref{claim(II)-1}, we deduce that the map
\[ \Ho^0(X, \mathcal{L} \otimes I_{Y_1}) \to \Ho^0(Y_2, (\mathcal{L} \otimes {\mathcal{N}_{Y_1/X}^*})|_{Y_2}) \]
is surjective. Now, in the following diagram,
\[ \begin{tikzcd}
0 \arrow{r} &\Ho^0(X, \mathcal{L} \otimes I_{Y_1}) \arrow{r} \arrow{d} & \Ho^0(X, \mathcal{L} \otimes I_{Y_2}) \arrow{d} \arrow{r} & \Ho^0(Y_1, \mathcal{L}|_{Y_1} \otimes I_{Y_2/Y_1}) \arrow{r} \arrow{d} & 0 \\
0 \arrow{r} & \Ho^0(Y_2, (\mathcal{L} \otimes {\mathcal{N}_{Y_1/X}^*})|_{Y_2}) \arrow{r} \arrow{d} & \Ho^0(Y_2, \mathcal{L}|_{Y_2} \otimes \mathcal{N}_{Y_2/X}^*) \arrow{r} & \Ho^0(Y_2, \mathcal{L}|_{Y_2} \otimes \mathcal{N}_{Y_2/Y_1}^*) \arrow{d} & \\
& 0 & & 0 &
\end{tikzcd} \]
the five lemma implies the surjectivity of the vertical morphism in the middle.
\end{pf}

\begin{claim} \rm
The map
\[ \Ho^0(X, p^*\stsh_{\mathbb{P}^s}(\alpha) \otimes \stsh_p(\beta) \otimes I_{Y_3}) \to \Ho^0(Y_3, q^*\stsh_{\mathbb{P}^s}(\alpha) \otimes \stsh_q(\beta) \otimes \mathcal{N}_{Y_3/X}). \]
is surjective.
\end{claim}

\begin{pf} \rm
By the same argument as in Claim \ref{claim(II)-2}, it is enough to show that
\begin{align*}
&\Ho^1(Y_2, \mathcal{L}|_{Y_2} \otimes I_{Y_3/Y_2}) = 0, \\
&\Ho^1(Y_2, \mathcal{L}|_{Y_2} \otimes I^2_{Y_3/Y_2}) = 0, \\
\text{and} ~&\Ho^1(Y_2, \mathcal{L}|_{Y_2} \otimes \mathcal{N}_{Y_2/X}^* \otimes I_{Y_3/Y_2}) = 0.
\end{align*}
We have $\mathcal{N}_{Y_2/X}^* = \bigoplus_{k=1,2} (p^*\stsh_{\mathbb{P}^s}(a_{\lambda_k}) \otimes \stsh_p(-1))|_{Y_2}$ and $I_{Y_3/Y_2} = p^*\stsh_{\mathbb{P}^s}(-1)|_{Y_2}$, and we obtain the vanishing of cohomology
\begin{align*}
\Ho^1(Y_2, \mathcal{L}|_{Y_2} \otimes I_{Y_3/Y_2}) &= \Ho^1(Y_2, (p^*\stsh_{\mathbb{P}^s}(\alpha-1) \otimes \stsh_p(\beta))|_{Y_2}) =0,
\end{align*}
since $\beta \geq 0$, $\alpha -1 \geq s- 2s' -1 \geq - s + 1$, and $Y_2$ is a $\mathbb{P}^{r-2}$-bundle over $\mathbb{P}^s$. 
Similarly, we have
\begin{align*}
\Ho^1(Y_2, \mathcal{L}|_{Y_2} \otimes I^2_{Y_3/Y_2}) &= \Ho^1(Y_2, (p^*\stsh_{\mathbb{P}^s}(\alpha-2) \otimes \stsh_p(\beta))|_{Y_2}) =0,
\end{align*}
and
\begin{align*}
\Ho^1(Y_2, \mathcal{L}|_{Y_2} \otimes \mathcal{N}_{Y_2/X}^* \otimes I_{Y_3/Y_2}) &= \bigoplus_{k=1,2} \Ho^1(Y_2, (p^*\stsh_{\mathbb{P}^s}(\alpha+a_{\lambda_k} -1) \otimes \stsh_p(\beta-1))|_{Y_2}) \\
&=0.
\end{align*}
Note that if $r=2$, then $\beta = 2$. \qed
\end{pf}
\end{pf}


\begin{thebibliography}{99}

\bibitem[Ba99]{Batyrev99} V. V. Batyrev, \textit{On the classification of toric Fano 4-folds},  J. Math. Sci. (New York), \textbf{94}(1) (1999) ,1021--1050.

\bibitem[Be79]{Beilinson79} A. Beilinson, \textit{Coherent sheaves on $\mathbb{P}^n$ and problems in linear algebra}, Funct. Anal. Appl. \textbf{12}(3) (1978), 68--69.

\bibitem[Bo90]{Bondal90} A. Bondal, \textit{Representation of associative algebras and coherent sheaves}, Math. USSR Izv. \textbf{34}(1) (1990), 23--42.

\bibitem[BK90]{BK90} A. Bondal, M. Kapranov, \textit{Representation functors, Serre functors, and reconstructions}, USSR Izv. \textbf{35}(3) (1990), 519--541.

\bibitem[BT09]{BT09} A. Bernardi, S. Tirabassi, \textit{Derived categories of toric Fano 3-folds via the Frobenius
morphism}, Matematiche (Catania), \textbf{64}(2) (2009), 117--154.

\bibitem[CM04]{CM04} L. Costa, R.M. Mir\'{o}-Roig, \textit{Tilting sheaves on toric varieties}, Math. Z. \textbf{248}(4) (2004), 849--865.

\bibitem[CM12]{CM12} L. Costa, R.M. Mir\'{o}-Roig, \textit{Derived categories of toric varieties with small Picard number}, Cent. Eur. J. Math. \textbf{10}(4) (2012), 1280--1291.

\bibitem[CLS]{CLS} D. Cox, J. Little, H. Schenck, \textit{Toric Varieties}, Graduate Studies in Mathematics, Vol. 124, Amer. Math. Soc., Providence, RI, 2011.

\bibitem[DLM09]{DLM09} A. Day, M. Laso\'{n}, M. Micha{\l}ek, \textit{Derived category of toric varieties with Picard number three}, Matematiche (Catania) \textbf{64}(2) (2009), 99--116.

\bibitem[Ef14]{Efimov14} A. Efimov, \textit{Maximal length of exceptional collections of line bundles}, J. London Math. Soc. \textbf{90}(2) (2014), 350--372.

\bibitem[HP06]{HP06} L. Hille, M. Perling, \textit{A counterexample to {K}ing's conjecture}, Compos. Math. \textbf{142}(6) (2006), 1507--1521.

\bibitem[Ka06]{Kawamata06} Y. Kawamata, \textit{Derived categories of toric varieties}, Michigan Math. J. \textbf{54}(3) (2006), 517--535.

\bibitem[Ki97]{King97} A. D. King, \textit{Tilting bundles on some rational surfaces}, preprint, 1997.

\bibitem[Kl88]{Kleinschmidt88} P. Kleinschmidt, \textit{A classification of toric varieties with few generators}, Aequationes Math. \textbf{35} (1988), 254--266.

\bibitem[LM11]{LM11} M. Laso\'{n}, M. Micha{\l}ek, \textit{On the full, strongly exceptional collections on toric varieties with Picard number three}, Collect. Math. \textbf{62} (2011), 275--296.

\bibitem[M11]{M11} M. Macha{\l}ek, \textit{Family of counterexamples to King's conjecture}, C. R. Math. Acad. Sci. Paris, \textbf{349} (2011), 67--69.

\bibitem[Or93]{Orlov93} D. Orlov, \textit{Projective bundles, monoidal transformations and derived categories of coherent sheaves}, Math. USSR Izv. \textbf{38} (1993), 133--144.

\bibitem[Pr15]{Prabhu-Naik15} N. Prabhu-Naik, \textit{Tilting bundles on toric Fano fourfolds}, J. Algebra, \textbf{471} (2017), 348--398.

\bibitem[Sa00]{Sato00} H. Sato, \textit{Toward the classification of higher-dimensional toric Fano varieties}, Tohoku Math. J. (2), \textbf{52}(3) (2000),383--413.

\bibitem[Ue14]{Uehara14} H. Uehara, \textit{Exceptional collections on toric Fano threefolds and birational geometry}, Internat. J. Math., \textbf{25}(7) (2014), 1450072, 32.

\end{thebibliography}
\end{document}